\documentclass{article}
\usepackage{indentfirst,latexsym,bm,amsmath,amsthm,amsfonts,array}
\begin{document}
\title{On complete Lie algebras}

\author{BinYong,Hsie
 \\{\small LMAM, Department of Mathematics, PeKing University,}
\\ {\small BeiJing, P.R.China 100871}
\\{\small e-mail: byhsie@math.pku.edu.cn}}
\date{}
\maketitle{}
\begin{abstract}
In this paper, the author gives two methods to construct complete
Lie algebras. Both methods show that the derivation algebras of
some Lie algebras are complete.
\end{abstract}

\footnotetext
{

MR(2000) Subject Classification: 17B45, 17B55, 17B56}

\section{\bf{Introduction}}
A finite dimensional Lie algebra $\mathfrak{g}$ is called a
complete Lie algebra if its center $C( \mathfrak{g} )$ is trivial
and all of its derivations are inner, i.e. $\mathrm{Der}(
\mathfrak{g} )=\mathrm{ad}( \mathfrak{g} )$.

There are many beautiful results on complete Lie algebras. A
famous result is E.V.Stitzinger's derivations towers theorem,
which says, for any finite dimensional Lie algebra $\mathfrak{h}$
with trivial center, there is a positive integer $n$ such that
$\mathrm{Der}^{n}(\mathfrak{h})$ is a complete Lie algebra. In
section 2, we construct some Lie algebras whose derivation
algebras are complete. We show that if $\mathfrak{g}$ is a Lie
algebra with $\mathfrak{b}$ a torus on it and $\tau$ the
commutator of $\mathfrak{b}$ in $\mathrm{Der}(\mathfrak{g})$ such
that $(\mathfrak{g},\mathfrak{b})$ is a non-degenerate pair, then
$\tau\times_{t}\mathfrak{g}$ is isomorphic to the derivation
algebra of $\mathfrak{b}\times_{t}\mathfrak{g}$ and is complete.
Especially, if $\mathfrak{b}$ is a maximal torus, then
$\mathfrak{b}\times_{t}\mathfrak{g}$ is complete. There are many
unsolved questions on complete Lie algebras. For example, one may
ask when a nilpotent Lie algebra is the nilpotent radical of some
solvable complete Lie algebra. We give part answer to this
question. In fact we show a nilpotent Lie algebra of
non-degenerate type is the nilpotent radical of some solvable
complete Lie algebra.

In section 3, at first, we give a method to calculate the
derivation algebras of the full graphs of Lie algebras of certain
type. Then we apply it to Heisenberg algebra to construct Lie
algebras which are not complete but whose derivation algebras are
complete.

\section{Complete Lie algebras associated to pairs of non-degenerate type}
\subsection{Notations and concepts}
Let $\mathfrak{g}$ be a Lie algebra, and let
$\mathrm{Der}(\mathfrak{g})$ denote its derivation algebra. A Lie
subalgebra $\mathfrak{b}$ of $\mathrm{Der}(\mathfrak{g})$ is
called a (maximal) torus on $\mathfrak{g}$ if $\mathfrak{b}$ is a
(maximal) commutative subalgebra whose elements are all
semisimple.

For a given torus $\mathfrak{b}$ on $\mathfrak{g}$ and any
$\alpha\in\mathfrak{b}^{*}$, put
$$\mathfrak{g} _{\alpha} = \{ g
\in \mathfrak{g}:\:D(g)=\alpha(D)g,\: \forall D \in \mathfrak{b}
\}.$$ Then $\mathfrak{g}$ has a decomposition
\begin{equation}
\mathfrak{g}=\bigoplus_{\alpha\in\Delta}
\mathfrak{g}_{\alpha}.\label{decomposition}
\end{equation}

A pair $(\mathfrak{g},\mathfrak{b})$ is called a non-degenerate
pair if all $\alpha$ in formula (\ref{decomposition}) are nonzero.

If $\mathfrak{g}$ is nilpotent and $\mathfrak{b}_{i}(i=1,2)$ are
two maximal toruses on $\mathfrak{g}$, then there is some
$\theta\in \mathrm{Aut} \mathfrak{g}$ such that
$\theta(\mathfrak{b}_{1})=\mathfrak{b}_{2}$. Hence
$(\mathfrak{g},\mathfrak{b}_{1})$ is a non-degenerate pair if and
only if $(\mathfrak{g},\mathfrak{b}_{2})$ is. Then we say that
$\mathfrak{g}$ is a nilpotent Lie algebra of non-degenerate type
or a non-degenerate nilpotent Lie algebra.

If there is a homomorphism $\varphi$ from Lie algebra
$\mathfrak{s}$ to $\mathrm{Der}(\mathfrak{g})$, we can define a
Lie algebra $\mathfrak{s}\times_{\varphi}\mathfrak{g}$ via
\begin{equation}
 [(s_{1},g_{1}),(s_{2},g_{2})]=([s_{1},s_{2}],s_{1}(g_{2})-s_{2}(g_{1})+[g_{1},g_{2}]).\label{demi-product}
\end{equation}
If $\mathfrak{s}=\mathrm{Der}(\mathfrak{g})$,
$\mathfrak{s}\times_{\varphi}\mathfrak{g}$ is called the full
graph of $\mathfrak{g}$ and is denoted by
$\mathfrak{f}(\mathfrak{g})$. If $\mathfrak{g}$ is center free,
then both $\mathrm{Der}(\mathfrak{g})$ and
$\mathfrak{f}(\mathfrak{g})$ are center free.

\subsection{Complete Lie algebras associated to pairs of non-degenerate type}
In this subsection, we assume $\mathfrak{b}$ is a torus on
$\mathfrak{g}$ such that $(\mathfrak{g},\mathfrak{b})$ is a
non-degenerate pair. Let $\tau$ be the commutator of
$\mathfrak{b}$ in $\mathrm{Der}(\mathfrak{g})$. Let
$\mathfrak{h}_{1}$ and $\mathfrak{h}$ denote $\mathfrak{b}
\times_{t} \mathfrak{g}$ and $\tau \times \mathfrak{g}$,
respectively. Then $\mathfrak{h}_{1}$ and $\mathfrak{h}$ have
decompositions
$$\mathfrak{h}_{1}=\mathfrak{b}\bigoplus \sum_{\alpha\in\Delta} \mathfrak{g}_{\alpha}$$
and
$$\mathfrak{h}=\tau\bigoplus \sum_{\alpha\in\Delta} \mathfrak{g}_{\alpha},$$
respectively.

We identity $\mathfrak{b}$ with $\{ (b,0):b\in\mathfrak{b} \}$ in
$\mathfrak{h}_{1}$, $\tau$ with $\{(t,0):t\in\tau\}$ in
$\mathfrak{h}$, and $\mathfrak{g}$ with
$\{(0,g):g\in\mathfrak{g}\}$ in $\mathfrak{h}_{1}$ and
$\mathfrak{h}$ respectively.

{\bf{Lemma 1}}. Let $\mathfrak{g}$, $\mathfrak{b}$, $\tau$ be as
above.

(i) If $D$ is a derivation of $\mathfrak{h}$ such that
$D(\mathfrak{b}) \subseteq \tau$, then $D=\mathrm{ad}(t)$ for some
$t\in\tau$.

(ii) If $D$ is a derivation of $\mathfrak{h}_{1}$ such that
$D(\mathfrak{b}) \subseteq \mathfrak{b}$, then $D|\mathfrak{b}=0$
and $D|\mathfrak{g}=t$ for some $t\in\tau$.

{\bf Proof.} We need only to prove (i), since (ii) can be proved
similarly.

We notice that each $\mathfrak{g}_{\alpha}$ is stable under
$\tau$, since $\tau$ commutes with $\mathfrak{b}$.

Assume $D e_{\alpha}=t_{1}+\sum_{\beta\in\Delta}x_{\beta}$
($t_{1}\in \tau,\: x_{\beta}\in\mathfrak{g}_{\beta}$) for a
nonzero element $e_{\alpha}$ in $\mathfrak{g}_{\alpha}$ with
$\alpha$ in $\Delta$. We have
\begin{equation}
\alpha(b)(t_{1}+\sum_{\beta\in\Delta}x_{\beta})=[Db,e_{\alpha}]+\sum_{\beta\in\Delta}\beta(b)x_{\beta}\;\;\;(b\in\mathfrak{b})
\end{equation}
since $D[b,e_{\alpha}]=[Db,e_{\alpha}]+[b,De_{\alpha}]$.
Therefore, we obtain $t_{1}=0$, $x_{\beta}=0\;(\beta\neq \alpha)$
and
\begin{equation}
[Db,e_{\alpha}]=0\;\;\forall\:\alpha\in\Delta,\;b\in\mathfrak{b}.
\label{1}
\end{equation}
Then we get $Db=0$ for any $b$ in $\mathfrak{b}$ and
$D(\mathfrak{g}_{\alpha})\subseteq\mathfrak{g}_{\alpha}$ for any
$\alpha$ in $\Delta$. Hence, $D|\mathfrak{g}$ commutates with
$\mathrm{ad}(\mathfrak{b})|\mathfrak{g}$ and so there exists a $t$
in $\tau$ such that
$$D|\mathfrak{g}=\mathrm{ad}(t)|\mathfrak{g}.$$
Let $D_{1}$ be $D-\mathrm{ad}(t)$, then
$D_{1}|\mathfrak{b}\bigoplus\mathfrak{g}=0$.

For any $t$ in $\tau$ and any $g$ in $\mathfrak{g}$, we have
$$[D_{1}t,g]=D_{1}[t,g]-[t,D_{1}g]=0.$$
Therefore, we obtain $D_{1}t=0$ as desired. \qed

{\bf Lemma 2.} Let
$\mathfrak{g},\mathfrak{b},\tau,\mathfrak{h}_{1}$ and
$\mathfrak{h}$ be as above.

(i). For any derivation $D$ of $\mathfrak{h}$, there exist
$x_{\alpha}\in\mathfrak{g}_{\alpha}$ ($\alpha\in\Delta$) such that
\begin{equation}
Db=b'+\sum_{\alpha\in\Delta}\alpha(b)x_{\alpha}\;\forall\:
b\in\mathfrak{b}, \end{equation} with some $b'\in\tau$.

(ii). For any derivation $D$ of $\mathfrak{h}_{1}$, there exist
$x_{\alpha}\in\mathfrak{g}_{\alpha}$ ($\alpha\in\Delta$) such that
\begin{equation}
Db=b'+\sum_{\alpha\in\Delta}\alpha(b)x_{\alpha}\forall \:b\in
\mathfrak{b}\end{equation} with some $b'\in \mathfrak{b}$.

{\bf Proof.} As in the proof of Lemma 1, we need only to prove
(i), since (ii) can be proved similarly.

For any $b$ in $\mathfrak{b}$ assume that
\begin{equation}
Db=b'+\sum_{\alpha\in\Delta}y_{\alpha}(b)
\end{equation}
with $b'\in \tau$ and $y_{\alpha}(b)\in\mathfrak{g}_{\alpha}$.
Then, for $ b_{1}$ in $\mathfrak{b}$, we have
\begin{equation}-\sum_{\alpha\in\Delta}\alpha(b_{1})y_{\alpha}(b)+\sum_{\alpha\in\Delta}\alpha(b)y_{\alpha}(b_{1})=0.\end{equation}
So
\begin{equation}\alpha(b_{1})y_{\alpha}(b)=\alpha(b)y_{\alpha}(b_{1})\;\forall \alpha\in \Delta, \;b,b_{1}\in \mathfrak{b}.\label{2}\end{equation}

If $\alpha(b)=0$, we find a $b_{1}$ in $\mathfrak{b}$ with
$\alpha(b_{1})\neq 0$, then from the above formula we obtain
$y_{\alpha}(b)=0$. If $y_{\alpha}(b)=0$ for all
$b\in\mathfrak{b}$, let $x_{\alpha}=0$. If there exists a
$b\in\mathfrak{b}$ with $y_{\alpha}(b)\neq 0$, then $\alpha(b)\neq
0$ and let $x_{\alpha}= \frac{1}{\alpha(b)} y_{\alpha}(b).$ Then
from formula (\ref{2}) we get
\begin{equation}y_{\alpha}(b_{1})=
\frac{\alpha(b_{1})}{\alpha(b)}y_{\alpha}(b)=\alpha(b_{1})x_{\alpha}\end{equation}
as desired.\qed

We can now prove the following theorem.

{\bf Theorem 1.} Let $\mathfrak{g}$ be a Lie algebra,
 and $\mathfrak{b}$ be a torus on $\mathfrak{g}$ such that
$(\mathfrak{g},\mathfrak{b})$ is a non-degenerate pair. Let $\tau$
be the commutator of $\mathfrak{b}$ in
$\mathrm{Der}(\mathfrak{g})$. Let $\mathfrak{h}$ and
$\mathfrak{h}_{1}$ denote $\tau\times_{t}\mathfrak{g}$ and
$\mathfrak{b}\times_{t}\mathfrak{g}$ respectively. Then
$\mathfrak{h}$ is equal to $\mathrm{Der}(\mathfrak{h}_{1})$ and is
a complete Lie algebra.

{\bf Proof.} (i). At first we prove
$\mathfrak{h}=\mathrm{Der}(\mathfrak{h}_{1})$.

Let $D$ be a derivation of $\mathfrak{h}_{1}$. From lemma 2 (ii)
we obtain
$$Db=b'+\sum_{\alpha\in\Delta} \alpha(b)x_{\alpha}\;\forall\: b\in \mathfrak{b}$$
with some $b'\in \mathfrak{b}$. Therefore
$D+\mathrm{ad}(\sum\limits_{\alpha\in\Delta}x_{\alpha})$ is a
derivation of $\mathfrak{h}_{1}$ such that
\begin{equation}
(D+\mathrm{ad}(\sum_{\alpha\in\Delta}x_{\alpha}))b=b'\in\mathfrak{b}\;\forall
\:b\in\mathfrak{b}.\end{equation} Lemma 1(ii)
shows\begin{equation}
D+\mathrm{ad}(\sum_{\alpha\in\Delta})x_{\alpha})|\mathfrak{b}=0\end{equation}
and
\begin{equation}
D+\mathrm{ad}(\sum_{\alpha\in\Delta}x_{\alpha})|\mathfrak{g}=t\label{4}
\end{equation}
for some $t\in\tau$.

For any $(t_{0},g_{0})$ in $\mathfrak{h}$, we define a derivation
$D_{(t_{0},g_{0})}$ of $\mathfrak{h}_{1}$ via
\begin{equation}
D_{(t_{0},g_{0})}((b,g))=(0,t_{0}(g)-b(g_{0})+[g_{0},g])\;\;\;\forall\:(b,g)\in\mathfrak{h}_{1}.\label{3}
\end{equation}
Then we get a homomorphism from $\mathfrak{h}$ to
$\mathrm{Der}(\mathfrak{h}_{1})$. It is easy to see
$D_{(t_{0},g_{0})}=0$ if and only if $(t_{0},g_{0})=0$. Thus we
can regard $\mathfrak{h}$ as a Lie subalgebra of
$\mathrm{Der}(\mathfrak{h}_{1})$ via (\ref{3}). From formula
(\ref{4}), we see that every derivation $D$ of $\mathfrak{h}_{1}$
lies in $\mathfrak{h}$. Therefore, $\mathfrak{h}$ is equal to
$\mathrm{Der}(\mathfrak{h}_{1})$.

(ii). We are going to prove $\mathfrak{h}$ is complete.

It is obvious that $C(\mathfrak{h})=0$. So it is sufficient to
prove $\mathrm{Der}(\mathfrak{h})=\mathrm{ad}(\mathfrak{h})$.

Let $D$ be a derivation of $\mathfrak{h}$. From lemma 2(i), we
obtain\begin{equation}
Db=b'+\sum_{\alpha\in\Delta}\alpha(b)x_{\alpha}\;\;\forall\: b\in
\mathfrak{b}
\end{equation}
with some $b'\in \tau$. Therefore, we have
\begin{equation}(D+\mathrm{ad}(\sum_{\alpha\in\Delta}x_{\alpha}))b=b'\in\tau\;\forall\:
b\in \mathfrak{b}.\end{equation} Then lemma 1(i) tells us that
there exists some $t\in\tau$ such that
$$D+\mathrm{ad}(\sum_{\alpha\in\Delta}x_{\alpha})=\mathrm{ad}(t)$$
So $D$ is in $\mathrm{ad}(\mathfrak{h})$.\qed

Now, we apply the above theorem to solvable complete Lie algebra.

A nilpotent Lie algebra is called completable if it is the
nilpotent radical of some solvable complete Lie algebra. We have
the following proposition.

{\bf Proposition 1.}  Let $\mathfrak{n}$ be a nilpotent Lie
algebra of non-degenerate type. Then $\mathfrak{n}$ is
completable.

{\bf Proof.} Let $\mathfrak{b}$ be a maximal torus over
$\mathfrak{n}$. Then its commutator $\tau$ in
$\mathrm{Der}(\mathfrak{n})$ is just $\mathfrak{b}$ itself.
Therefore, $\mathfrak{b}\times_{t}\mathfrak{n}$ is a solvable
complete Lie algebra.\qed

We use proposition 1 to construct some examples of completable
nilpotent lie algebras. Let $\mathfrak{g}$ be a Lie algebra and
let $n$ be a positive integer. We will define a Lie algebra
$\mathfrak{g}^{\oplus n}=\{g=(g(1),...,g(n))|\:g(1),...,g(n)\:
\mathrm{are\:elements \: of\:\mathfrak{g}}\}$ via
$$[g,g'](k)=\sum_{i+j=k}[g(i),g'(j)],\;\;\forall\:g,g'\in \mathfrak{g}^{\oplus n},\; 1\leq i,j,k\leq
n.$$ We show $\mathfrak{g}^{\oplus n}$ is a non-degenerate
nilpotent Lie algebra. Let $D$ be a derivation over
$\mathfrak{g}^{\oplus n}$ such that $D(g)(i)=i g(i)$ for any $g\in
\mathfrak{g}^{\oplus n}$, then $(\{tD| t \in
\mathbb{R}\},\mathfrak{g}^{\oplus n})$ is a non-degenerate pair.
Hence, $\mathfrak{g}^{\oplus n}$ is non-degenerate and so it is
completable.

\section{Some relation between complete Lie algebra and Heisenberg algebra}
\subsection{A general result}
In this section, we are going to study full graphs of Lie
algebras. Let $\mathfrak{g}$ be a Lie algebra and $\mathcal{S}$ be
a complete Lie algebra with a homomorphism $\varphi$ from
$\mathcal{S}$ to $\mathrm{Der}(\mathfrak{g})$. Let
 $\mathfrak{h}$ denote $\mathcal{S}\times_{\varphi}\mathfrak{g}$.

We define $Z_{s}(\mathfrak{g})$ via \begin{equation}
Z_{s}(\mathfrak{g})=\{g\in C(\mathfrak{g}):s(g)=0\forall s\in
\mathcal{S}\}. \end{equation}

 By a simple calculate we obtain the following lemma.

{\bf Lemma 3.} The center of $\mathfrak{h}$ is $\{(0,g):g\in
Z_{s}(\mathfrak{g})\}$.

Let $F_{s}(\mathfrak{g})$ be the set $\{D\in
\mathrm{Der}(\mathfrak{g}):[D,D_{1}]\in
\mathrm{ad}(\mathfrak{g}),\;\forall\: D_{1}\in
\varphi(\mathcal{S})\}$ which is a Lie subalgebra of
$\mathrm{Der}(\mathfrak{g})$. It is obvious that
$\mathrm{ad}(\mathfrak{g})$ is its ideal. We have a natural
induced homomorphism $\tilde{\varphi}$ from $\mathcal{S}$ to
$\mathrm{Der}(\mathrm{Der}(\mathfrak{g}))$. If $\mathcal{S}$ is
just $\mathrm{Der}(\mathfrak{g})$, let $F(\mathfrak{g})$ denote
$F_{s}(\mathfrak{g})$. Let $\widetilde{\mathfrak{h}}$ denote
$\mathcal{S}\times_{\tilde{\varphi}}F_{s}(\mathfrak{g})$.

{\bf Lemma 4.} If $C(\mathfrak{g})=0$, then
$\mathrm{Der}(\mathfrak{h}) \supseteq \widetilde{\mathfrak{h}}$.

Furthermore, if one of the following conditions (a), (b) and (c)
is satisfied, then
$\mathrm{Der}(\mathfrak{h})=\widetilde{\mathfrak{h}}$.

(a). $\mathrm{Der}(\mathfrak{g})\supseteq \mathcal{S}\supseteq
\mathrm{ad}(\mathfrak{g})$,

(b). $[\mathfrak{g},\mathfrak{g}]=\mathfrak{g}$,

(c). $\mathfrak{g}=\mathfrak{b}\times \mathfrak{n}$ is a solvable
Lie algebra with $(\mathfrak{n},\mathfrak{b})$ a non-degenerate
pair, and $\mathcal{S}$ annihilates $\mathfrak{b}$.

{\bf Proof.} For any inner derivation $D=\mathrm{ad}(g)$ of
 $\mathfrak{g}$, let $I(D)$ denote $g$.
For any $(s,D)$ in $\widetilde{\mathfrak{h}}$, a derivation
$\mathrm{D}_{(s,D)}$ on $\mathfrak{h}$ is defined via
\begin{equation}
\mathrm{D}_{(s,D)}(s_{1},g)=([s,s_{1}],s(g)+D(g)+I([D,\varphi(s_{1})])\;\;\forall\:(s_{1},g)\in\mathfrak{h}.\end{equation}
By a simple calculate, we get the following assertions

 $\;\;(\mathrm{i}).\;\;\;\;\mathrm{D}_{(s,D)}$ is a derivation,

 $\;\;(\mathrm{ii}).\;\;\;
[\mathrm{D}_{(s_{1},D_{1})},\mathrm{D}_{(s_{2},D_{2})}]=\mathrm{D}_{[(s_{1},D_{1})(s_{2},D_{2})],}$\\
and

 $\;\;(\mathrm{iii}).\:\; \mathrm{D}_{(s,D)}=0$ if and only if $s=0$ and
 $D=0$.\\
Therefore, $\widetilde{\mathfrak{h}}$ is a Lie subalgebra of
$\mathrm{Der}(\mathfrak{h})$.

Now we begin to prove
$\mathrm{Der}(\mathfrak{h})=\widetilde{\mathfrak{h}}$ if at least
one of the conditions (a), (b) and (c) is satisfied. Let
$\mathrm{D}_{1}$ be a derivation of $\mathfrak{h}$. Write
$\mathrm{D}_{1}=(T'_{1},T_{2})$, where $T'_{1}$ takes values in
$\mathcal{S}$ and $T_{2}$ takes values in $\mathfrak{g}$. It is
obvious that $T'_{1}|\mathcal{S}$ is a derivation on
$\mathcal{S}$. Since $\mathcal{S}$ is complete, we see
$T'_{1}|\mathcal{S}=\mathrm{ad}(s)|\mathcal{S}$ with some $s\in
\mathcal{S}$. Let $\mathrm{D}=\mathrm{D}_{1}-\mathrm{D}_{(s,0)}$,
then $\mathrm{D}=(T_{1},T_{2})$ with $T_{1}|\mathcal{S}=0$. We are
going to show $T_{1}=0$. By
$$\mathrm{D}([(0,g_{1})(s,g_{2})])=[\mathrm{D}(0,g_{1}),(s,g_{2})]+[(0,g_{1}),\mathrm{D}(s,g_{2})],$$
we obtain
\begin{equation}
-T_{1}(s(g_{1}))+T_{1}([g_{1},g_{2}])=[T_{1}(g_{1}),s]\;\;\;\forall\:
s\in \mathcal{S}, g_{1},g_{2}\in
\mathfrak{g}.\label{5}
\end{equation}
We see \begin{equation}
T_{1}|[\mathfrak{g},\mathfrak{g}]=0,\label{6}\end{equation} and
\begin{equation} -T_{1}(s(g))=[T_{1}(g),s].\label{7}\end{equation}

If (b) holds, it is obvious $T_{1}=0$. Next we assume (a) holds.
From formulas (\ref{6}) and (\ref{7}), we obtain \begin{equation}
[T_{1}(g),s]=0\;\;\;\forall\: s\in \mathrm{ad}(\mathfrak{g}),g\in
\mathfrak{g},\end{equation} since $\mathrm{ad}(\mathfrak{g})$ is
contained in $\mathcal{S}$. Hence, $T_{1}$ annihilates
$\mathfrak{g}$. Now we assume (c) holds. Condition (c) implies
that $\mathfrak{g}=\mathfrak{b}+[\mathfrak{g},\mathfrak{g}]$. From
formula (\ref{7}), we get\begin{equation}
[T_{1}(b),s]=0\;\;\forall\;b\in \mathfrak{b}\;,s\in\mathcal{S}.
\end{equation}
Since $C(\mathcal{S})=0$, we see that $T_{1}$ annihilates
$\mathfrak{b}$. Since $T_{1}|[\mathfrak{g},\mathfrak{g}]=0$, the
proof of $T_{1}=0$ is completed.

We have already obtained $\mathrm{D}=(0,T_{2})$. The following
equalities hold\begin{equation}
T_{2}([(s_{1},0),(s_{2},0)])=[(s_{1},0),T_{2}(s_{2},0)]+[T_{2}(s_{1},0),(s_{2},0)]\label{40.1}\end{equation}
\begin{equation}
T_{2}([(s,0),(0,x)])=[(s,0),T_{2}(0,x)]+[T_{2}(s,0),(0,x)]\label{40.2}\end{equation}
\begin{equation}
T_{2}([(0,x_{1}),(0,x_{2})])=[(0,x_{1}),T_{2}(0,x_{2})]+[T_{2}(0,x_{1}),(0,x_{2})].\label{40.3}\end{equation}
 From formula (\ref{40.3}), we see $ T_{2}|\mathfrak{g}\in
\mathrm{Der}(\mathfrak{g})$ and we may assume $
T_{2}|\mathfrak{g}=D.$ From (\ref{40.2}) we obtain that
\begin{equation}
\mathrm{ad}(T_{2}(s))=[D,\varphi(s)]\;\forall\;
s\in\mathcal{S}.\label{42.1}\end{equation} Therefore $D$ belongs
to $F_{s}(\mathfrak{g})$ as desired.
 \qed

From lemma 3 and lemma 4, it is easy to obtain the following
theorem.

{\bf Theorem 2.} Let $\mathfrak{g}$ be a Lie algebra with trivial
center, and assume that $\mathrm{Der}(\mathfrak{g})$ is a complete
lie algebra. Then the derivation algebra of the full graph
$\mathfrak{f}(\mathfrak{g})$ of $\mathfrak{g}$ is
$\mathrm{Der}(\mathfrak{g})\times_{\varphi} F(\mathfrak{g})$.
Therefore, $\mathfrak{f}(\mathfrak{g})$ is a complete Lie algebra
if and only if $F(\mathfrak{g})=\mathfrak{g}$.

For any integer $n$, let $\mathfrak{f}^{n}(\mathfrak{g})$ denote
$\mathfrak{f}(\mathfrak{f}^{n-1}(\mathfrak{g}))$ by induction.

The following corollary can be easily proved by induction.

{\bf Corollary.} Let $\mathfrak{g}$ be a non-complete Lie algebra
with trivial center, and assume that $\mathrm{Der}(\mathfrak{g})$
is a complete Lie algebra such that
$[\mathrm{Der}(\mathfrak{g}),\mathrm{Der}(\mathfrak{g})]\subseteq
\mathrm{ad}(\mathfrak{g})$. Then we have the following assertions

\ \ (i).
$\mathrm{Der}(\mathfrak{f}^{n}(\mathfrak{g}))=\mathfrak{f}^{n}(\mathrm{Der}(\mathfrak{g}))$
is complete but $\mathfrak{f}^{n}(\mathfrak{g})$ is not,

\ (ii). The center of $\mathfrak{f}^{n}(\mathfrak{g})$ is trivial,

(iii). $[\mathrm{Der}(\mathfrak{f}^{n}(\mathfrak{g})),
\mathrm{Der}(\mathfrak{f}^{n}(\mathfrak{g}))]\subseteq
\mathrm{ad}(\mathfrak{f}^{n}(\mathfrak{g}))$,

(iv).
$\mathrm{dim}(\mathrm{Der}(\mathfrak{f}^{n}(\mathfrak{g})))-\mathrm{dim}(\mathfrak{f}^{n}(\mathfrak{g}))=\mathrm{dim}(\mathrm{Der}(\mathfrak{g}))-\mathrm{dim}(\mathfrak{g})$.


\subsection{Apply the general result to Heisenberg algebras}

In this section, let $\mathfrak{g}$ be a Heisenberg algebra such
that $[\mathfrak{g},\mathfrak{g}]=C(\mathfrak{g})$ and
$\mathrm{dim}(C(\mathfrak{g}))=1$.

Let $c$ be a nonzero element of $C(\mathfrak{g})$. Then there is
an anti-symmetric quadratic form $\psi$ on $\mathfrak{g}$ of rank
$\mathrm{dim}(\mathfrak{g})-1$ such that for any $s,\:y$ in
$\mathfrak{g}$
\begin{equation} [x,y]=\psi(x,y)c.
\end{equation}

Then dimension of $(\mathfrak{g})$ is equal to $2N+1$ for some
positive integer $N$. We may choose a basis
$\{x_{1},...,x_{N};y_{1},...,y_{N};c\}$ such that for any $1\leq
i,j\leq N$,\begin{equation}
\psi(x_{i},x_{j})=\psi(y_{i},y_{j})=\psi(c,x)=0,\end{equation} and
\begin{equation} \psi(x_{i},y_{j})=\delta_{ij}.\end{equation}

In \cite{BZ}, the following three propositions are stated.

{\bf{Proposition 2}}. $\mathrm{Der}(\mathfrak{g})$ is simple
complete Lie algebra, and has the following decomposition
$$\mathrm{Der}(\mathfrak{g})=\mathcal{S}\times_{t}(\mathfrak{b}\times_{t}\mathfrak{n}),$$
where, $\mathcal{S}\cong
\mathrm{sp}(2n,\mathbb{C}),\;\mathfrak{n}=\mathrm{ad}(\mathfrak{g})\cong\mathfrak{g}/\mathbb{C}c$
is a irreducible $\mathcal{S}$-module and
$\mathrm{dim}(\mathfrak{b})=1$.

{\bf{Proposition 3}}. The full graph $\mathfrak{f}(\mathfrak{g})$
of $\mathfrak{g}$ is not a complete Lie algebra, but its center is
trivial.

{\bf{Proposition 4}}. The derivation algebra
$\mathrm{Der}(\mathfrak{f}(\mathfrak{g}))$ of
$\mathfrak{f}(\mathfrak{g})$ is a complete Lie algebra with
$\mathrm{dim}(\mathrm{Der}(\mathfrak{f}(\mathfrak{g})))=\mathrm{dim}(\mathfrak{f}(\mathfrak{g}))+1$.

From proposition 3 and proposition 4, we obtain
$$[\mathrm{Der}(\mathfrak{f}(\mathfrak{g})),\mathrm{Der}(\mathfrak{f}(\mathfrak{g}))]\subseteq \mathrm{ad}(\mathfrak{f}(\mathfrak{g})).$$

By corollary of theorem 2, we get

{\bf Theorem 3.} Let $\mathfrak{g}$ be a Heisenberg algebra, then
for any positive integer $n$,

(i). The center of $\mathfrak{f}^{n}(\mathfrak{g})$ is trivial,

(ii).
$\mathrm{Der}(\mathfrak{f}^{n}(\mathfrak{g}))=\mathfrak{f}^{n-1}(\mathrm{Der}(\mathfrak{f}(\mathfrak{g})))$
is complete but $\mathfrak{f}^{n}(\mathfrak{g})$ is not complete,

(iii).
$[\mathrm{Der}(\mathfrak{f}^{n}(\mathfrak{g})),\mathrm{Der}(\mathfrak{f}^{n}(\mathfrak{g}))]\subseteq
\mathrm{ad}(\mathfrak{f}^{n}(\mathfrak{g})) $,

(iv).
$\mathrm{dim}(\mathrm{Der}(\mathfrak{f}^{n}(\mathfrak{g})))-\mathrm{dim}(\mathfrak{f}^{n}(\mathfrak{g}))=1$.

\end{document}